\documentclass[oneside,a4paper, 12pt]{article}
\usepackage{anysize}
\usepackage{retraction}
\usepackage{biber}

\hypersetup{pdftitle={Short retractions of CAT(1) spaces},
pdfauthor={Alexander Lytchak and Anton Petrunin}}

\begin{document}
\title{Short retractions of CAT(1) spaces}
\author{Alexander Lytchak and Anton Petrunin}
\date{}
\maketitle

\begin{abstract}
We construct short retractions of a CAT(1) space to its small convex subsets.
This construction provides an alternative geometric description of an analytic tool introduced by Wilfrid Kendall.

Our construction uses a \emph{tractrix flow} which can be defined as a gradient flow for a family of functions of certain type.
In an appendix we prove a general existence result for gradient flows of time-dependent locally Lipschitz semiconcave functions, which is of independent interest.
\end{abstract}

\section{Introduction}

Recall that a  subset $K$ in a metric space ${\spc{U}}$ is called \emph{weakly convex} if any two points $x,y\in K$ can be connected by a minimizing geodesic in $K$.

Let ${\spc{U}}$ be a metric space and $K\subset {\spc{U}}$.
A distance nonexpanding map $f\:{\spc{U}}\to K$ such that $f(x)=x$ for any $x\in K$ is called a \emph{short retraction to $K$}.
If in addition a local Lipschitz constant of $f$ is strictly less than 1 at any point $x\notin K$, 
then we say that $f$ is a \emph{strictly short retraction from ${\spc{U}}$ to $K$}.

\begin{thm}{Theorem}\label{thm:retraction:Phi}
Let $\spc{U}$ be a complete length $\CAT(\kappa)$ space.
Suppose $K$ is a weakly convex closed subset in $\spc{U}$ and there is $p\in \spc{U}$ such that $|p-x|\le \tfrac\pi2$ for any point $x\in K$.

\begin{enumerate}[(a)]
 \item If $p\in K$ and $\kappa\le 1$, then there is a short retraction 
$\Phi\:\spc{U}\to K$.
\item If $\kappa<1$, then there is a strictly short retraction 
$\Phi\:\spc{U}\to K$.
\end{enumerate}
\end{thm}

This statement is a generalization of the following well known statement about $\CAT(0)$ spaces:
\emph{If ${\spc{U}}$ is a complete length $\CAT(0)$ space and $K$ is a closed convex subset in ${\spc{U}}$,
then the closest point projection ${\spc{U}}\to K$ is a short retraction.
Moreover, if $\spc{U}$ is a $\CAT(\kappa)$ space for some $\kappa<0$, then  the closest point projection is a strictly short retraction}.

The theorem and a small trick imply the following:

\begin{thm}{Corollary}\label{cor}
Let $\spc{U}$ be a complete length $\CAT(\kappa)$ space.
Denote by $\Delta$ the diagonal in $\spc{U}\times \spc{U}$;
that is, $\Delta=\{\,(x,x)\in \spc{U}\times \spc{U}\,\}$.

Suppose there is a point $p\in \spc{U}$ such that $|p-x|\le \tfrac\pi2$ for any point $x\in \spc{U}$.
\begin{enumerate}[(a)]
\item
If $\kappa\le 1$, then there is a short retraction $\Psi\:\spc{U}\times \spc{U}\to \Delta$.
\item If $\kappa<1$, then there is a strictly short retraction $\Psi\:\spc{U}\times \spc{U}\to \Delta$.
\end{enumerate}

\end{thm}

It is well known that if $\spc{U}$ is a complete length $\CAT(0)$ space, then the midpoint map $\spc{U}\times \spc{U}\to \spc{U}$ is $\tfrac1{\sqrt{2}}$-Lipschitz and therefore it induces a short retraction $\spc{U}\times \spc{U}\to\Delta$. 
The corollary provides an analogous statement for $\CAT(1)$ spaces.

\parbf{Motivation.}
In \cite[(4.1)]{kendall}, Wilfrid Kendall observed that if $\spc{B}$ is a regular geodesic ball of radius $r<\tfrac\pi2$ in a manifold with sectional curvature at most 1, then, for an appropriate choice of constant $\lambda$, the function
\[(x,y)\mapsto 
\frac{1+\lambda-\cos|x-y|_{\spc{B}}}{\cos|p-x|_{\spc{B}}\cdot \cos|p-y|_{\spc{B}}}
\]
has convex level sets in ${\spc{B}}\times {\spc{B}}$.
He also shows the existence of a nonnegative convex function on ${\spc{B}}\times {\spc{B}}$ that vanishes only on the diagonal \cite[(4.2)]{kendall}.
These observations became a useful tool to study the Dirichlet problem and its relatives;
they allowed to extend a number of results from Hadamard manifolds to Riemannian manifolds of small size
and more generally to $\CAT(1)$ spaces~\cite{yokota,BFHMSZ,fuglede,serbinowski,lytchak-stadler}. 

Our original goal was to make this tool transparent for geometers.
Corollary~\ref{cor} can be considered as a more geometric version of this tool.
While Kendall's condition is optimal for uniqueness and regularity questions, the existence statements can be derived from Theorem~\ref{thm:retraction:Phi} in a
slightly greater generality, as  we are going to explain now.

We will need the following definition, introduced by Stefan Wenger \cite{Wenger-1comp};
for the definitions of ultrafilters and ultracompletions we refer to \cite{Wenger-1comp,guo-wenger,akp}.

A metric space ${\spc{U}}$ is called \emph{$1$-complemented} if for some \emph{non-principal ultrafilter} $\omega$  there exists a short retraction of the ultracompletion ${\spc{U}}^{\omega}$ to~${\spc{U}}$.
Examples of $1$-complemented spaces include all proper spaces, all CAT(0) spaces and all $L^p$ spaces for $1\leq p\leq \infty$ \cite[Proposition 2.1]{guo-wenger}.

Recall that if ${\spc{U}}$ is $\CAT(\kappa)$, then so is~${\spc{U}}^\omega$. 
Applying these observations together with Theorem~\ref{thm:retraction:Phi}, we obtain

\begin{thm}{Theorem}\label{thm:complemented}
Let ${\spc{U}}$ be a complete length $\CAT(1)$ space.
Assume  there exists some $p\in {\spc{U}}$ such that $|p-x|\le \tfrac\pi2$ for any point $x\in {\spc{U}}$.
Then ${\spc{U}}$ is $1$-complemented.
\end{thm}

Let us list a few existence results which follow from the theorem, assuming that the space ${\spc{U}}$ is as above:
\begin{enumerate}[(a)]
\item\label{dirichlet}   The existence of a solution $u$ of Dirichlet problem on the minimization of energy 
in $W^{1,2} (\Omega, {\spc{U}})$ on any Lipschitz domain $\Omega$ in a Riemannian manifold with prescribed trace $tr(u)$; see \cite{KS} and \cite[Theorem 1.4]{guo-wenger}.
\item\label{current} The existence of a minimal integral $k$-current filling any prescribed boundary in ${\spc{U}}$; 
see \cite{Ambrosio} and \cite[Theorem 3.3]{Wenger-1comp}.
\item   The existence of a conformally parametrized disc $u:D\to {\spc{U}}$ of minimal area for a given boundary curve $\gamma$, which is a Jordan curve of finite length in ${\spc{U}}$;
see \cite{LWplateau} and \cite[Theorem 1.2]{guo-wenger}.
\item\label{center} For any Radon measure $\mu$ on ${\spc{U}}$ there exists a center of mass $x\in {\spc{U}}$ for the measure~$\mu$ \cite{Sturm, yokota}.
\end{enumerate}

If in the theorem we assume strict inequality $|p-x|< \tfrac\pi2$, then the existence results are known  in all the cases (\ref{dirichlet}--\ref{center}).
Moreover, the uniqueness holds true under this stronger assumption in  the cases (\ref{dirichlet}), (\ref{current}), and (\ref{center}); see
\cite{yokota,serbinowski}.
In our boundary case uniqueness definitely fails;
for example geodesics between points in a round hemisphere are not unique.

The uniqueness in each case can be shown using Corollary \ref{cor}.
Indeed if there are different solutions of one of these problems, then their product in ${\spc{U}}\times {\spc{U}}$ does not lie in the diagonal.
The latter contradicts the existence of the strictly short retraction $\Psi$ provided by Corollary \ref{cor}.

\parbf{About the proofs.}
We use a new tool which we call $r$-\emph{tractrix flow}, a special  time-dependent gradient flow.
It gives a family of maps $\phi_t$ for a given rectifiable curve $t\z\mapsto\gamma(t)$.
The important properties of the tractrix flow are collected in Proposition~\ref{prop-def}.
In particular, (1) if $\spc{U}$ is $\CAT(1)$ and $r\le \tfrac\pi2$, then $\phi_t$ is short for any $t$, 
and (2) if $r< \tfrac\pi2$, then the local Lipschitz constant of $\phi_t$ at $p$ is strictly less than 1 if $p\ne \phi_t(p)$.

In the proof of Theorem~\ref{thm:retraction:Phi}, the tractrix flow is applied in a space obtained by gluing to $\spc{U}$ a spherical cone over $K$;
this space is $\CAT(1)$ by Reshetnyak's gluing theorem.
In Appendix~\ref{Another way} we indicate another way of proving Theorem~\ref{thm:retraction:Phi}.

In the proof of Corollary~\ref{cor} the additional trick consists in identifying the product space $\spc{U}\z\times \spc{U}$ with a subset of the spherical join $\spc{U}\star \spc{U}$ and applying Theorem~\ref{thm:retraction:Phi} to the latter.

The tricks in both proofs show that it is useful to consider singular spaces even in the case when the original space $\spc{U}$ is smooth;
this is a powerful freedom of Alexandrov's world.
More involved examples of such arguments are given by Dmitry Burago, Sergei Ferleger, and Alexey Kanonenko \cite{BFK}, Paul Creutz~\cite{creutz}, and Stephan Stadler~\cite{stadler}. 

\parbf{Acknowledgements.}
We thank Christine Brenier,
Nicola Gigli,
and the anonymous referee for helpful comments.
A. Lytchak was partially supported by DFG grant SPP 2026.
A. Petrunin was partially supported by Simons foundation collaboration grant for mathematicians 584781.

\section{Tractrix flow}\label{sec:Tractrix flow}

For $\CAT(\kappa)$ spaces, we will follow the conventions in \cite{akp}.

First let us describe the tractrix flow informally.
Suppose that two points $p$ and $q$ in ${\spc{U}}$ are connected to each other by a thread of fixed length $r$.
Imagine that the point $q$ follows the curve $\gamma$ and drags $p$ if the thread is tight; 
if the thread is not tight, then $p$ does not move.
Then the trajectory of the point $p$ will be called $r$-\emph{tractrix of $p$ with respect to~$\gamma$}.
The family of maps $\phi_t$ that send the initial position of $p$ to its position at the time $t$ will be called the $r$-\emph{tractrix flow} defined by $\gamma$.

More formally, suppose $\gamma\:[a,b]\to \spc{U}$ is a 1-Lipschitz curve. 
An $r$-tractrix with respect to $\gamma$ is defined as a gradient curve for the time-dependent family of functions 
\[f_t=- \max\{r,\dist_{\gamma(t)}\};\]
here $\dist_{x}$ denotes the distance function from the point $x$.
We also assume that the initial point lies in $\bar B(\gamma(a),r)$.
(We denote by $\bar B(x,r)$ and $B(x,r)$ the closed and open balls of radius $r$ centered at $x$.) 

The $r$-\emph{tractrix flow} with respect to $\gamma$ is defined as a family of maps
\[\phi_t\:\bar B(\gamma(a),r)\z\to\bar B(\gamma(t),r)\]
whose trajectory  $t\mapsto \phi_t(p)$ is the $r$-tractrix starting at $p\in \bar B(\gamma(a),r)$.

The following proposition includes the properties of the tractrix flow which will be used further in the paper.

\begin{thm}{Proposition}\label{prop-def}
Let $\gamma\:[a,b]\to \spc{U}$ be a 1-Lipschitz curve in a complete length $\CAT(\kappa)$ space $\spc{U}$ for some $\kappa\le 1$
and $r<\pi$.
Set $\bar B_t=\bar B(\gamma(t),r)$.
Then the $r$-tractrix flow $\phi_t\:\bar B_a\to\bar B_t$ is uniquely defined.
Moreover
\begin{enumerate}[(a)]
\item  \label{approx} For any $t$, the map $\phi_t$ is a limit of compositions $\theta_{t_n} \circ \dots\circ \theta_{t_0}$
for $\delta \to 0$, where $a=t_0<\dots<t_n=t$ is any partition of $[a,t]$ with $|t_i-t_{i-1}|<\delta$
and where $\theta_{t_i} : \bar B_{t_{i-1}}\to  \bar B_{t_i}$  denotes the closest point projection.
\item \label{sharafutdinov} If the family of balls $\bar B_t$ is decreasing in $t$ (that is, if $\bar B_{t_1}\supset \bar B_{t_2}$ for $t_1<t_2$), then $\phi_b$ is a strong deformation retraction from $\bar B_a$ to $\bar B_b$.
\item \label{non-strict} If $r=\tfrac\pi2$, then $\phi_t$ is short for any $t$;
 \item\label{strict} If $r=\tfrac\pi2$ and  $\kappa<1$, then there is a positive constant $\eps$ such that the local Lipschitz constant of $\phi_t$ at $p$ is bounded above by $\exp(-\eps\cdot\ell)$, 
 where $\ell=|p-\phi_t(p)|_{\spc{U}}$.
\end{enumerate}
\end{thm}

Historically the first relative of the tractrix flow
is the so called \emph{Sharafutdinov's retraction} \cite{sharafutdinov} --- a family of maps associated to a continuous family of convex sets (in our case these sets are the balls $\bar B_t$). 
Second relative is the \emph{pursuer flow} introduced and studied by Stephanie Alexander, Richard Bishop, Robert Ghrist and Chanyoung Jun \cite{ABG,jun-thesis,jun,jun:grad}.

Time-dependent gradient flows were studied by Chanyoung Jun \cite{jun-thesis,jun:grad}, by Lucas C. F. Ferreira and Julio C. Valencia-Guevara \cite{ferreira-valencia}, and by Alexander Mielke, Riccarda Rossi, and Giuseppe Savar\'{e} \cite{mielke-rossi-savare}.
Unfortunately Proposition~\ref{prop-def} does not follow directly from the results in these papers; for this reason we provide in an appendix a short proof of the existence of
gradient flows of Lipschitz time-dependent family of semiconcave functions in $\CAT(\kappa)$ spaces.

\parit{Proof.}
Consider $f_t=- \max\{r,\dist_{\gamma(t)}\}$ as a family of functions defined in $B(\gamma(t),r\z+\delta)$ for sufficiently small $\delta>0$.
Note that the family $f_t$ is Lipschitz.
By $\CAT(1)$-comparison, each $f_t$ is $\lambda$-concave for a fixed $\lambda$.
Moreover if $r<\tfrac\pi2$, then $\lambda=0$ and if $r=\tfrac\pi2$, then $\lambda\to 0$ as $\delta\to 0$. 

Consider the map $\phi_t\:\alpha(a)\mapsto \alpha(t)$, where $\alpha$ is a $f_t$-gradient curve.
By \ref{prop:time-dependent}, if $\phi_t(p)$ is defined, then it is unique.

Consider the function $\ell(t)\z=|\phi_t(p)-\gamma(t)|$.
By the definition of the flow, we have that $\ell'\le 0$ if $\ell> r$.
It follows that $\phi_t$ is defined for all $t$ and maps $\bar B_a$ to $\bar B_t$.

\parit{(\ref{approx}).} 
Given a partition $a=t_0<t_1<\dots<t_n=t$ with $|t_i -t_{i-1}|<\delta$, consider a locally constant approximation $\hat f_t$ of the family $f_t$ defined by $\hat f_t=f_{t_i}$ if $t_i\le t  < t_{i+1}$.
Denote by $\hat\phi_t$ the corresponding flow.

Given $p\in \bar B_a$, set $p_i=\hat\phi_{t_i}(p)$.
Observe that $p_{i}=\theta_{t_i}(p_{i-1})$ for each $i$.

By the distance estimate (\ref{Distance estimate}) the flow $\hat\phi_t$ converges to $\phi_t$ as the partition gets finer and finer, hence the result.

\parit{(\ref{sharafutdinov}).} It is sufficient to notice that $\phi_t(p)=p$ if
 $|p-\gamma(s)| \leq r$ for all $a\leq s \leq t$.

\parit{(\ref{non-strict}).} 
Applying the distance estimate (\ref{Distance estimate}) for $s=0$, we get that 
\[|\phi_t(p)-\phi_t(q)|\le |p-q|\cdot e^{\lambda\cdot (t-a)}\]
for any $p,q\in \bar B_a$, $t\ge a$.
If $r\le \tfrac\pi2$, then the inequality holds for arbitrary $\lambda> 0$;
hence (\ref{non-strict}) follows.

\parit{(\ref{strict}).} 
 The proof of the strict inequality follows directly from  (\ref{approx})
 and the following general consequence of the $\CAT(\kappa)$ comparison:

 There exists some $\eps>0$ such that 
 the closest point projection $\theta:\bar B(w,\frac \pi 2 +\eps)\z\to\bar B(w,\frac \pi  2)$  in any $\CAT(\kappa)$ is strictly short and satisfies 
 \[|\theta (p) -\theta (q)| \leq e^{-\eps\cdot |\theta (p) -p|}\cdot|p-q|.\]
\qedsf

\section{Proofs}\label{sec:proofs}

Recall that spherical join $\spc{U}\star\spc{V}$ of two metric spaces $\spc{U}$ and $\spc{V}$
is defined as the unit sphere (equipped with the angle metric) in the product of Euclidean cones $\Cone \spc{U}\times \Cone\spc{V}$. 
If diameters of $\spc{U}$ and $\spc{V}$ do not exceed $\pi$, then $\spc{U}\star\spc{V}$
can be defined as a metric space that admits an onto map $\iota\:\spc{U}\times\spc{V}\times[0,\tfrac\pi2]\to \spc{U}\star\spc{V}$ such that
\[
\begin{aligned}
|\iota(u_1,v_1,t_1)&-\iota(u_2,v_2,t_2)|_{\spc{U}\star\spc{V}}=
\\
&=\arccos[\sin t_1\cdot\sin t_2\cdot \cos|u_1-u_2|_{\spc{U}}+\cos t_1\cdot \cos t_2\cdot \cos|v_1-v_2|_{\spc{V}}].
\end{aligned}
\eqlbl{join-formula}
\]

Recall that the join of two $\CAT(1)$ spaces is $\CAT(1)$ \cite[Corollary 3.14]{bridson-haefliger}.

\parit{Proof of \ref{thm:retraction:Phi}.}
Consider the join of  $K$ with a one-point space, $\spc{J}=K\star \{s\}$.
Since $\spc{J}$ is a $\CAT(1)$ space,
by Reshetnyak's gluing theorem \cite[8.9.1]{akp}, the space $\spc{W}$ glued from ${\spc{U}}$ and $\spc{J}$ along $K$ is a $\CAT(1)$ space;
moreover ${\spc{U}}$ and $\spc{J}$ are convex subsets in $\spc{W}$.

\begin{wrapfigure}{o}{50 mm}
\vskip-0mm
\centering
\includegraphics{mppics/pic-2}
\end{wrapfigure}

Let $\gamma$ be the geodesic in $\spc{W}$ from $p$ to the pole $s$ of $\spc{J}$.
Set $\bar B_t=\bar B(\gamma(t),\tfrac\pi2)_{\spc{W}}$, then
\begin{itemize}
\item $\bar B_0=\bar B(p,\tfrac\pi2)_{\spc{U}}\cup \spc{J}$,
\item $\bar B_{\frac\pi2}=\spc{J}$, in particular $\bar B_{\frac\pi2}\cap \spc{U}=K$,
\item the family $\bar B_t$ is decreasing in $t$.
\end{itemize}
To see the last statement, note that $\spc{J}\subset \bar B_0$.
Further, by definition of a spherical join we have $|x-p|_{\spc{W}} \le |x\z-\gamma (s)|_{\spc{W}}$ for any $x\in K$.
By construction of the metric on $\spc{W}$ the same inequality holds for any
$x\in \spc{U}$,
Therefore $\bar B_s\subset \bar B_0$.
Finally, by $\CAT(1)$ comparison, $\bar B_t\supset \bar B_s\cap \bar B_0=\bar B_s$ for any $s\ge t \ge 0$.

According to \ref{prop-def}(\ref{sharafutdinov}), the $\tfrac\pi2$-tractrix flow $\phi_t$ is a strong deformation retraction of $\bar B_0$ to $\bar B_{\frac\pi2}$.
By \ref{prop-def}(\ref{non-strict}) $\phi_{\frac\pi2}$ is a short.
If $\kappa<1$, then by \ref{prop-def}(\ref{strict}), $\phi_{\frac\pi2}$ is a strictly short retraction.

Since ${\spc{U}}$ is $\CAT(1)$,
given a point $x\in B(p,\pi)_{\spc{U}}$ there is unique geodesic $\gamma_x$ parametrized by its length from $p$ to $x$. 
By $\CAT(1)$ comparison, the map 
\[\Theta(x)=
\begin{cases}
p&\text{if\ }|p-x|_{\spc{U}}\ge \pi,
\\
\gamma_x(\pi-|p-x|_{\spc{U}})&\text{if\ }|p-x|_{\spc{U}}< \pi.
\end{cases}
\]
is a short retraction of ${\spc{U}}$ to $\bar B(p,\tfrac\pi2)_{\spc{U}}=\bar B_0\cap \spc{U}$.
Moreover $\Theta$ is strictly short retraction if $\kappa<1$.

Therefore the composition $\Phi=\phi_{\frac\pi2}\circ\Theta$ induces a short retraction of ${\spc{U}}$ to $K$
which is strictly short if $\kappa<1$.

Finally, we need to take care of the case $\kappa<1$ and $p\notin K$.
Denote by $\bar p\in K$ the closest point to $p$; by $\CAT(\kappa)$ comparison it exists and unique.
Note that $|\bar p-x|_{\spc{U}}< |p-x|_{\spc{U}}$ for any $x\in K$;
therefore $K\subset B(\bar p,\tfrac\pi2)$. 
It remains to apply the construction above with $\bar p$ instead $p$.
\qeds

\parit{Proof of \ref{cor}.}
Consider the spherical join $\spc{U}\star\spc{U}$ and the map $\iota$ described at the beginning of the section. 
Note that \ref{join-formula} implies that the map $(u,v)\mapsto \iota(u,v,\tfrac\pi4)$
induces a length preserving map 
\[\Theta\:{\tfrac1{\sqrt{2}}}\cdot(\spc{U}\times\spc{U})\to\spc{U}\star\spc{U}.\]
In particular, $\Theta$ is short.

Note that the diagonal $\tfrac1{\sqrt{2}}\cdot\Delta$ is a convex set in $\tfrac1{\sqrt{2}}\cdot(\spc{U}\times\spc{U})$.
Moreover \ref{join-formula} implies that the restriction of $\Theta$ to $\tfrac1{\sqrt{2}}\cdot\Delta$ is distance preserving.
In particular, the image $K=\Theta(\tfrac1{\sqrt{2}}\cdot\Delta)$ is a weakly convex set in $\spc{U}\star\spc{U}$.

Further note that $|q-y|_{\spc{U}\star\spc{U}}\le \tfrac\pi2$ for any $y\in \spc{U}\star\spc{U}$ and $q=\Theta(p,p)$.
Applying \ref{thm:retraction:Phi}, we get a short retraction $\Phi\:\spc{U}\star\spc{U}\to K$.
Since $\Theta$ is short, it induces the needed short retraction $\Psi\:\spc{U}\times\spc{U}\to \Delta$.

Finally, by \ref{thm:retraction:Phi}, if $\kappa<1$, then $\Phi$ is a strictly short retraction and therefore so is $\Psi$.
\qeds

\appendix
\section{Time-dependent gradient flow}

Here we prove the existence, uniqueness and contractivity  of the gradient flow for a time-dependent family of functions.
The proof relies mostly on the corresponding statements  for time-independent families --- with minor conditions on the space.
The same  proof works nearly without changes for spaces with lower curvature bound
and it should work in nearly any space with well defined angles between geodesics starting at one point.

\parbf{Time-independent flow.}
Suppose $\spc{U}$ is a complete length $\CAT(\kappa)$ space.

For a locally Lipschitz semiconcave function $f$ defined on an open set $\Dom f\subset \spc{U}$, the differential $d_pf\:\T_p\to\RR$ is defined at each point $p\in \Dom f$;
it is a concave, Lipschitz, and positive-homogeneous of degree 1 function on the tangent space at $p$. 

Further, for any point $p\in \Dom f$ there is unique tangent vector $u\in \T_p$
such that the following two conditions
\[
\begin{aligned}
\langle u,w\rangle &\ge d_{p}f(w),
\\
\langle u,u\rangle &= d_{p}f(u)
\end{aligned}
\eqlbl{<u,w>}
\]
hold for any tangent vector $w\in \T_p$.%
\footnote{Here \emph{tangent vector} means \emph{element of tangent cone}, we use this term despite the tangent cone is not a vector space.
The scalar product $\langle u,w\rangle$ is defined as $|u|\cdot|w|\cdot\cos\theta$ where $\theta$ is the angle between the vectors.}
The vector $u$ is called the \emph{gradient of $f$ at $p$}; briefly $u=\nabla_pf$.

A locally Lipschitz map $t\mapsto \alpha(t)$ is called \emph{$f$-gradient curve} if it satisfies the following equation
\[\alpha^+(t)=\nabla_{\alpha(t)}f\eqlbl{grad-flow}\]
for any $t$. 
Here $\alpha^+(t)\in \T_{\alpha(t)}$ denotes the \emph{right velocity vector}; that is,
\[\alpha^+(t)=\lim_{\eps\to 0+} \frac{\log_{\alpha(t)}[\alpha(t+\eps)]}\eps,\]
where $v=\log_pq$ if the vector $v\in\T_p$ points form $p$ in the direction of $q$ and $|v|=|p-q|$.

The following proposition can be extracted from \cite[Theorem 1.7]{lytchak-open-map} or \cite{ohta-palfia}.

\begin{thm}{Proposition}\label{prop:time-independent}
Let ${\spc{U}}$ be a complete length $\CAT(\kappa)$ space and
$f$ a Lipschitz and semiconcave function defined on an open set $\Dom f\subset {\spc{U}}$.
Then for any point $p\in \Dom f$  
there is a unique $f$-gradient curve $t\mapsto\alpha(t)$ with initial point $\alpha(0)=p$ defined on a maximal semiopen interval $[0,T)$ for some $T\in(0,\infty]$.
Moreover if $T<\infty$ then $\alpha$ escapes from any closed set $K\subset \Dom f$.
\end{thm}

\parbf{Time-dependent gradient flow.}
Our next aim is to define a gradient flow for a time-dependent family of functions and prove an analog of Proposition~\ref{prop:time-independent} for this flow.
Let $f_t$ be a family of functions defined on open subsets $\Dom f_t$ of~$\spc{U}$.
More precisely, we assume that the parameter $t$ lies in a real interval $\II$ and 
\[\Omega=\set{(x,t)\in\spc{U}\times \II}{x\in\Dom f_t}\]
is an open subset in $\spc{U}\times \II$.

A family of functions $f_t$ is called \emph{Lipschitz} if 
the function $(x,t)\mapsto f_t(x)$ is $L$-Lipschitz for some constant $L$.

A family of functions $f_t$ will be called \emph{semiconcave} if 
the function $x\mapsto f_t(x)$ is $\lambda$-concave for each $t$.
A family $f_t$ is called \emph{locally semiconcave} if for each $(p_0,t_0)\in \Omega$ there is a neighborhood $\Omega'$ and $\lambda\in\RR$ such that the restriction of $f_t$ to $\Omega'$ is $\lambda$-concave. 

Note that one cannot expect that a direct generalization of equation \ref{grad-flow} holds for any family of functions $f_t$.

For example, consider a $1$-Lipschitz curve $\alpha$ in the real line. 
It is reasonable to assume that $\alpha$ is an $f_t$-gradient curve for the family $f_t(x)=-|x-\alpha(t)|$.
(Indeed $\alpha$ can be realized as a limit of  gradient curves for a family of functions obtained by smoothing $f_t$.)
On the other hand, $\alpha^+(t)$ might be undefined,
but even if it is defined, $\alpha^+(t)\ne0$ in general, while $\nabla_{\alpha(t)} f_t\equiv0$.

Instead we define \emph{$f_t$-gradient curve} as a Lipschitz curve $\alpha$ that satisfies the following inequality
for any point $p$, time $t$, and
small $\eps >0$:  
\[\dist_p\circ\alpha(t+\eps)\le \dist_p\circ\alpha(t)-\eps\cdot d_{\alpha(t)}f_t(\dir{\alpha(t)}p)+o(\eps),\eqlbl{def:tdflow}\]
where $\dir qp\in \T_{q}$ denotes a unit tangent vector at $q$ in the direction of $p$
(if there is no geodesic $[\alpha(t)\,p]$ then we impose no condition).

If $\alpha^+(t)=\nabla_{\alpha(t)}f_t$ for all $t$, then \ref{def:tdflow} holds;
it follows from \ref{<u,w>}.
On the other hand, the example above shows that the converse does not hold;
that is, \ref{def:tdflow} generalizes the definition~\ref{grad-flow}.
Our defining  inequality \ref{def:tdflow} is closely related to the so called \emph{evolution variational inequality} \cite[Thm 4.0.4(iii)]{ambrosio-gigli-savare}.

\begin{thm}{Distance estimate}\label{Distance estimate}
Let $f_t$ and $h_t$ be two families of $\lambda$-concave functions on a $\CAT(\kappa)$ space $\spc{U}$ and $s\ge 0$.
Assume $f_t$ and $h_t$ have common domain $\Omega\subset {\spc{U}}\times \RR$ and $|f_t(x)-h_t(x)|\le s$ for any $(x,t)\in \Omega$.
Assume $t\mapsto \alpha(t)$ and $t\mapsto \beta(t)$ are $f_t$- and $h_t$-gradient curves respectively defined on a common interval $t\in [a,b)$; set $\ell(t)\z=|\alpha(t)\z-\beta(t)|_{\spc{U}}$.
If for all $t$ a minimizing geodesic $[\alpha(t)\,\beta(t)]$ lies in $\set{x\in {\spc{U}}}{(x,t)\in \Omega}$, then
\[\ell'(t)\le \lambda\cdot\ell(t)+2\cdot s/\ell(t),
\eqlbl{eq:ell'=<}\]
assuming that the left hand side is defined.
Moreover
\[\ell(t)^2+\tfrac{2\cdot s}\lambda\le(\ell(a)^2+\tfrac{2\cdot s}\lambda)\cdot e^{2\cdot\lambda\cdot (t-a)}.
\eqlbl{eq:ell2=<}\]

In particular the inequalities hold for any $t\in\II$ if $\Omega\supset B(p,2\cdot r)\times \II$ and $\alpha(t),\beta(t)\z\in B(p, r)$ for any $t\in \II$.
\end{thm}

Note that if $f_t=h_t$ then $s=0$;
in this case the second inequality can be written as
\[\ell(t)\le \ell(a)\cdot e^{\lambda\cdot (t-a)}.\eqlbl{dist-est-s=0}\]
In particular it implies uniqueness of the future of gradient curve with given initial data.
This inequality also makes it possible to estimate the distance between two gradient curves for close functions.
In particular, it implies convergence of $f_t^n$-gradient curves if a sequence of $L$-Lipschitz and $\lambda$-concave families $f^n_t$ converges uniformly as $n\to \infty$. 

\parit{Proof.}
Fix a time moment $t$ and set $f=f_t$ and $h=h_t$.
Let $p$ be the midpoint of the geodesic $[\alpha(t)\beta(t)]$.
Let $\gamma\:[0,\ell(t)]\to \spc{U}$ be an arc length parametrization of $[\alpha(t)\beta(t)]$ from $\alpha(t)$ to $\beta(t)$.
Note that $d_{\alpha(t)}f(\dir{\alpha(t)}{p})$ is the right derivative of $f\circ\gamma$ at $0$
and $-d_{\alpha(t)}h(\dir{\beta(t)}p)$ is the left derivative of $h\circ\gamma$ at $\ell(t)$.
Since $f$ and $h$ are $\lambda$-concave,
\begin{align*}
f(\beta(t))&\le f(\alpha(t))+\ell(t)\cdot d_{\alpha(t)}f(\dir{\alpha(t)}{p}) +\tfrac12\cdot\lambda\cdot\ell(t)^2,
\\
h(\alpha(t))&\le h(\beta(t))+\ell(t)\cdot d_{\alpha(t)}h(\dir{\beta(t)}p) +\tfrac12\cdot\lambda\cdot\ell(t)^2,
\end{align*}
Adding these inequalities up and taking into account that $|f(x)-h(x)|<s$ for any $x$, we conclude that 
\[d_{\alpha(t)}f(\dir{\alpha(t)}{p})+d_{\alpha(t)}h(\dir{\beta(t)}p)\ge \lambda\cdot \ell(t)+2\cdot s/\ell(t).\]

Applying the triangle inequality and the definition of gradient curve at $p$, we get that
\begin{align*}
\ell(t+\eps)&=|\alpha(t+\eps)-\beta(t+\eps)|\le
\\
&\le |\alpha(t+\eps)-p|+|\beta(t+\eps)-p|\le 
\\
&\le |\alpha(t)-p|-\eps\cdot d_{\alpha(t)}f(\dir{\alpha(t)}{p})+|\beta(t+\eps)-p|-\eps\cdot d_{\beta(t)}h(\dir{\beta(t)}p)+o(\eps)=
\\
&=\ell(t)-\eps\cdot(\lambda\cdot \ell(t)+2\cdot s/\ell(t))+o(\eps)
\end{align*}
for $\eps>0$; hence the statement.

Since $\alpha$ and $\beta$ are Lipschitz, $t\mapsto \ell(t)$ is a Lipschitz function.
By Rademacher's theorem, its derivative $\ell'$ is defined almost everywhere and it satisfies the fundamental theorem of calculus.
Therefore \ref{eq:ell'=<} implies \ref{eq:ell2=<}.
\qeds

\begin{thm}{Proposition}\label{prop:def-time-dependent}
Let $\spc{U}$ be a complete length $\CAT(\kappa)$ space.
Suppose $f_t$ is a family of $\lambda$-concave functions for $t\in [a,b)$ and $\Dom f_t\supset B(z,2\cdot r)$ for some fixed $z\in\spc{U}$, $r>0$ and any~$t$.

Let $\alpha\:[a,b)\to B(z,r)$ be Lipschitz.
Then $\alpha$ is an $f_t$-gradient curve if and only if 
\[\dist_p\circ\alpha(t+\eps)\le \dist_p\circ\alpha(t)-\eps\cdot \left[\frac{f_t(p)-f_t\circ\alpha(t)}{|p-\alpha(t)|}-\tfrac\lambda2\cdot |p-\alpha(t)|\right]+o(\eps)\eqlbl{def:tdflow-plus}\]
for any $t\in [a,b)$ and $p\in B(z,r)\setminus \{\alpha (t)\}$.
\end{thm}

\parit{Proof.}
Note that geodesics $[\alpha(t)p]$ lie in $\Dom f_t$ for any $t$.

Since $f_t$ is $\lambda$-concave, we have 
\[d_{\alpha(t)}f_t(\dir{\alpha(t)}p)
\ge
\frac{f(p)-f\circ\alpha(t)}{|p-\alpha(t)|}-\tfrac\lambda2\cdot |p-\alpha(t)|.\]
Hence the only-if part follows.

Given a point $p\in \spc{U}$ and $t$,
consider a point $\bar p\in [\alpha(t)p]$.
Applying \ref{def:tdflow-plus} for $\bar p$ and the triangle inequality, we get
\[\dist_p\circ\alpha(t+\eps)
\le
\dist_p\circ\alpha(t)-\eps\cdot \left[\frac{f(\bar p)-f\circ\alpha(t)}{|\bar p-\alpha(t)|}-\tfrac\lambda2\cdot |\bar p-\alpha(t)|\right]+o(\eps).\]
By taking $\bar p$ close to $\alpha(t)$,
the value $\tfrac{f(\bar p)-f\circ\alpha(t)}{|\bar p-\alpha(t)|}-\tfrac\lambda2\cdot |\bar p-\alpha(t)|$ can be made arbitrary close to $d_{\alpha(t)}f_t(\dir{\alpha(t)}p)$.
Therefore, given $\delta>0$, the following inequality
\[\dist_p\circ\alpha(t+\eps)\le \dist_p\circ\alpha(t)-\eps\cdot d_{\alpha(t)}f_t(\dir{\alpha(t)}p)+\eps\cdot\delta\]
holds for all sufficiently small positive values $\eps$.
Therefore \ref{def:tdflow} holds.
\qeds

Now we are ready to formulate and prove an analog of Proposition~\ref{prop:time-independent} for time-dependent family.

\begin{thm}{Theorem}\label{prop:time-dependent}
Let $\spc{U}$ be a complete length $\CAT(\kappa)$ space and
$\{f_t\}$ a family of function defined on an open set $\Omega=\set{(x,t)\in \spc{U}\times\RR}{x\in \Dom f_t}$.
Suppose that $f_t$ is Lipschitz and locally semiconcave.
Then for any time moment $a$ and initial point $p\in \Dom f_a$ there is a unique $f_t$-gradient curve $t\mapsto\alpha(t)$ defined on a maximal semiopen interval $[a,b)$. 
Moreover if $b<\infty$, then $(\alpha(t),t)$ escapes from any closed set $K\subset \Omega$.
\end{thm}

\parit{Proof.}
Let $L$ be a Lipschitz constant of $f_t$.
Fix $b>a$ sufficiently small so that $\Dom f_t\z\supset B(p,\eps\cdot L)$ for any $t\in[a,b)$.
Consider sequence  $a=t_0<t_1\dots<t_n\z=b$ and a piecewise constant family of functions on $B(p,\eps\cdot L)$ defined by $\hat f_t=f_{t_i}$ if $t_i\le t<t_{i+1}$.

Note that $\hat f_t$ is time-independent on each interval $[t_i,t_{i+1})$.
Applying \ref{prop:time-independent} recursively on each interval $[t_i,t_{i+1})$, we get that the proposition holds for $\hat f_t$.
That is, there is a unique $\hat f_t$-gradient curve $\hat \alpha$ that starts at $p$ and defined on the interval $[a,b)$.

The distance estimates (\ref{Distance estimate}) show that as the partition gets finer, the gradient curves $\hat\alpha$ form a Cauchy sequence; denote its limit by $\alpha$.
Then
\begin{align*}
\dist_p\circ\hat\alpha(t+\eps)
&\le 
\dist_p\circ\hat\alpha(t)
-\eps\cdot \left[\frac{\hat f_t(p)-\hat f_t\circ\hat\alpha(t)}{|p-\alpha(t)|}-\tfrac\lambda2\cdot |p-\hat\alpha(t)|\right]
+o(\eps)
\\
&\le 
\dist_p\circ\hat\alpha(t)
-\eps\cdot \left[\frac{f_t(p)-f_t\circ\hat\alpha(t)-2\cdot s}{|p-\alpha(t)|}-\tfrac\lambda2\cdot |p-\hat\alpha(t)|\right]
+o(\eps)
\end{align*}
where 
\[s=\sup_{t,x} \{|f_t(x)-\hat f_t(x)|\}.\]
Since $s\to 0$ as $\hat\alpha\to \alpha$, we get that \ref{def:tdflow-plus} holds for $\alpha$;
that is, $\alpha$ is an $f_t$-gradient curve.

It proves a short time existence.
Applying this argument recursively we can find a gradient curve defined on a maximal interval $[a,b)$.
Uniqueness of this curve follows from the distance estimate \ref{dist-est-s=0}. 

Note that $\alpha$ is $L$-Lipschitz.
In particular, if $b<\infty$, then $\alpha(t)\to p'$ as $t\to b$.
If $(p',b)\in \Omega$, then we can repeat the procedure, otherwise $\alpha$ escapes from any closed set in $\Omega$. 
\qeds

\section{Another way}\label{Another way}

Here we indicate an alternative way to prove Theorem~\ref{thm:retraction:Phi}.
The idea is taken from the proof of Kirszbraun's theorem; see \cite[5.1]{akp-kirszbraun} or \cite[9.4.1 and 9.6.5]{akp}. 

\begin{wrapfigure}{o}{55 mm}
\centering
\includegraphics{mppics/pic-1}
\end{wrapfigure}

\parit{Proof of \ref{thm:retraction:Phi}.}
Set $\mathring{\spc{U}}=\Cone \spc{U}$.
Denote by $\mathring{K}$ the subcone of $\mathring{\spc{U}}$ spanned by $K$.
The space $\spc{U}$ is the unit sphere in $\mathring{\spc{U}}$ with angle metric;
that is, $\spc{U}=\set{z\in \mathring{\spc{U}}}{|z|=1}$ and $|u-v|_{\spc{U}}=\measuredangle[o\,{}^u_v]$, where $o$ denotes the tip of the cone $\mathring{\spc{U}}$.

Note that $\mathring{\spc{U}}$ is a $\CAT(0)$ space and $\mathring{K}$ forms a convex closed subset of $\mathring{\spc{U}}$.
In particular, for any point $x$ there is unique point $\hat x\in \mathring{K}$
that minimize the distance to $x$;
that is, $|\hat x-x|=\dist_K(x)$.

Consider the ray $\alpha_o\:t\mapsto t\cdot p$ in  $\mathring{\spc{U}}$.
Note that given $s\in \mathring{\spc{U}}$
the geodesics $[s\ t\cdot p]$ converge as $t\to\infty$ to a ray, 
say $\alpha_s\:[0,\infty)\to \mathring{\spc{U}}$ that is parallel to $\alpha_o$.

Note that if $|x|=1$, then $|\hat x|\le 1$.
By assumption for any $k\in \mathring{K}$ the function $t\mapsto |\alpha_k(t)|$ is monotonically increasing.
Therefore there is unique value $t_x\ge 0$ such that
$|\alpha_{\hat x}(t_x)|=1$;
set $x'=\alpha_{\hat x}(t_x)$ and define the map $\Phi$ as $x\mapsto x'$.

Evidently $x\mapsto x'$ is a retraction of $\spc{U}$ to $K$;
that is,
$\Phi(x)\in K$ for any $x\in \spc{U}$
and 
$\Phi(k)=k$ for any $k\in K$.

It remains to show that $x\mapsto x'$ is short which can be done by straightforward lengthy computations.
We leave it as an advanced exercise to the reader.
\qeds

{\small\sloppy

\printbibliography[heading=bibintoc]

}

\end{document}